\documentclass[11pt]{article}
\usepackage{amsmath}
\usepackage{amssymb}
\usepackage{amsthm}
\usepackage[usenames]{color}
\usepackage{amscd}

\usepackage[colorlinks=true,
linkcolor=webgreen,
filecolor=webbrown,
citecolor=webgreen]{hyperref}

\definecolor{webgreen}{rgb}{0,.5,0}
\definecolor{webbrown}{rgb}{.6,0,0}

\def\N{{\Bbb N}}
\def\Z{{\Bbb Z}}
\def\C{{\Bbb C}}
\def\1{{\bf 1}}

\def\ba{\begin{aligned}} \def\ea{\end{aligned}}
\def\br#1{\left(#1\right)}
\def\be{\begin{equation}} \def\ee{\end{equation}}
\def\c{(2\gamma-1)}
\def\ep{\varepsilon}
\def\h{\textstyle{\frac{1}{2}}}
\def\inv#1{\frac{1}{#1}}
\def\Int{\int\limits}
\def\d{\,\mathrm{d}}
\def\div{\,|\,}
\def\abs#1{\left|#1\right|}
\def\bo{\beta_0(K)}
\def\bi{\beta_1(K)}

\def\rank{\operatorname{rank}}

\def\lcm{\operatorname{lcm}}

\newtheorem{theorem}{Theorem}

\newtheorem{remark}{Remark}
\newtheorem{lemma}{Lemma}

\begin{document}

\title{\bf On the average number of subgroups of the group $\Z_m \times \Z_n$}
\author{Werner Georg Nowak and L\'aszl\'o T\'oth}
\date{}
\maketitle

\begin{abstract} Let $\Z_m$ be the group of residue classes modulo $m$. Let
$s(m,n)$ and $c(m,n)$ denote the total number of subgroups of the
group $\Z_m \times \Z_n$ and the number of its cyclic subgroups,
respectively, where $m$ and $n$ are arbitrary positive integers. We
derive asymptotic formulas for the sums $\sum_{m,n\le x} s(m,n)$,
$\sum_{m,n\le x} c(m,n)$ and for the corresponding sums restricted
to $\gcd(m,n)>1$, i.e., concerning the groups $\Z_m \times \Z_n$
having rank two.
\end{abstract}

{\sl 2010 Mathematics Subject Classification}: 11A25, 11N37, 20K01,
20K27

{\sl Key Words and Phrases}: cyclic group, finite Abelian group of rank
two, number of subgroups, number of cyclic subgroups, Dirichlet series,
asymptotic formula, Dirichlet's divisor problem


\section{Introduction}

Throughout the paper we use the notations: $\N:=\{1,2,\ldots\}$,
$\N_0:=\{0,1,2,\ldots\}$, $\Z_m$ is the additive group of residue
classes modulo $m$, $\phi$ is Euler's totient function, $\tau(n)$ is
the number of divisors of $n$, $\mu$ denotes the M\"obius function,
$\psi$ is the Dedekind function given by $\psi(n)=n\prod_{p\div n}
(1+1/p)$, $*$ stands for the Dirichlet convolution of arithmetic
functions, $\zeta$ is the Riemann zeta-function. Let $n=\prod_p p^{\nu_p(n)}$ be the prime
power factorization of $n\in \N$, where the product is over the primes $p$ and all but a finite number
of the exponents $\nu_p(n)$ are zero. Furthermore, let $\gamma_k$ ($k\in
\N_0$) denote the Stieltjes constants defined by
\begin{equation*}
\gamma_k :=\lim_{x\to \infty} \left(\sum_{n\le x}
\frac{(\log n)^k}{n} - \frac{(\log x)^{k+1}}{k+1} \right) \,,
\end{equation*}
where $\gamma_0=\gamma$ is the Euler-Mascheroni constant. We note
that the constants $\gamma_k$ are connected to the coefficients of
the Laurent series expansion of the function $\zeta(s)$ about its
pole $s=1$, namely,
\begin{equation*}
\zeta(s)= \frac1{s-1} +\sum_{k=0}^{\infty}
\frac{(-1)^k\gamma_k}{k!}(s-1)^k \,,
\end{equation*}
see, e.g., A.~Ivi\'c \cite[Th.\ 1.3]{Ivi1985}.

Consider the group $G :=\Z_m \times \Z_n$, where $m,n\in \N$ are
arbitrary. Note that $G$ is isomorphic to $\Z_{\gcd(m,n)}\times
\Z_{\lcm(m,n)}$. If $\gcd(m,n)=1$, then $G$ is cyclic, isomorphic to
$\Z_{mn}$. If $\gcd(m,n)>1$, then $G$ has rank two. We recall that a
finite Abelian group has rank $r$ if it is isomorphic to $\Z_{n_1}
\times \cdots \times \Z_{n_r}$, where $n_1,\ldots,n_r \in
\N\setminus \{1\}$ and $n_j \div n_{j+1}$ ($1\le j\le
r-1$). 
Let $s(m,n)$ and $c(m,n)$ denote the
total number of subgroups of the group $G$ and the number of its
cyclic subgroups, respectively.

Concerning general properties of the subgroup lattice of finite
Abelian groups, see R.~Schmidt \cite{Sch1994}, M.~Suzuki
\cite{Suz1951}. For every $m,n\in \N$ one has
\begin{equation} \label{s_decomp}
s(m,n)= \prod_p s(p^{\nu_p(m)},p^{\nu_p(n)}) \,,
\end{equation}
\begin{equation} \label{c_decomp}
c(m,n)= \prod_p c(p^{\nu_p(m)},p^{\nu_p(n)}) \,.
\end{equation}

For the $p$-group $\Z_{p^a} \times \Z_{p^b}$ of rank two, with $1\le
a\le b$, the following explicit formulas hold:
\begin{equation} \label{s_prime_pow}
s(p^a,p^b)=
\frac{(b-a+1)p^{a+2}-(b-a-1)p^{a+1}-(a+b+3)p+(a+b+1)}{(p-1)^2} \,,
\end{equation}
\begin{equation} \label{c_prime_pow}
c(p^a,p^b)= 2(1+p+p^2+\ldots +p^{a-1})+(b-a+1)p^a \,.
\end{equation}

The formula \eqref{s_prime_pow} was deduced, applying Goursat's
lemma for groups, by G.~C\u{a}lug\u{a}reanu \cite[Sect.\ 4]{Cal2004}
and J.~Petrillo \cite[Prop.\ 2]{Pet2011}, and using the concept of
the fundamental group lattice by M.~T\u{a}rn\u{a}uceanu \cite[Prop.\
2.9]{Tar2007}, \cite[Th.\ 3.3]{Tar2010}. Formula \eqref{c_prime_pow}
was given in \cite[Th.\ 4.2]{Tar2010}. Therefore, $s(m,n)$ and
$c(m,n)$ can be computed using \eqref{s_decomp}, \eqref{s_prime_pow}
and \eqref{c_decomp}, \eqref{c_prime_pow}, respectively. The
following more compact formulas were derived in \cite{HHTW2012} by a
simple elementary method: For every $m,n\in \N$,
\begin{equation} \label{eval_s_m_n}
\ba s(m,n) &=  \sum_{d \div m, e\div n} \gcd(d,e) \\ &= \sum_{d \div
\gcd(m,n)} \phi(d)\tau(m/d)\tau(n/d) \ea
\end{equation}
and
\begin{equation} \label{eval_c_m_n}
\ba c(m,n) &=  \sum_{d \div m, e\div n} \phi(\gcd(d,e)) \\ &= \sum_{d \div
\gcd(m,n)} (\mu*\phi)(d) \tau(m/d)\tau(n/d). \ea
\end{equation}

See also \cite{Tot2011,Tot2012} for a general identity concerning the
number of cyclic subgroups of an arbitrary finite Abelian group.

The identities \eqref{s_decomp} and \eqref{c_decomp} tell us that
the functions $(m,n)\mapsto s(m,n)$ and $(m,n)\mapsto c(m,n)$ are
multiplicative, viewed as arithmetic functions of two variables.
This property follows also from the first formulas of
\eqref{eval_s_m_n} and \eqref{eval_c_m_n}, respectively. See
\cite[Sect.\ 2]{HHTW2012}. Therefore, $n\mapsto s(n):=s(n,n)$
(sequence \cite[A060724]{OEIS}) and $n\mapsto c(n):=c(n,n)$
(sequence \cite[A060648]{OEIS}) are multiplicative functions of a
single variable and for every $n\in \N$ one has
\begin{equation} \label{sss}
s(n) =\sum_{d\div n} \tau(d) \psi(n/d) \,,
\end{equation}
\begin{equation} \label{ccc}
c(n)=\sum_{d\div n} \psi(d) \,,
\end{equation}
see \cite[Sect.\ 3]{HHTW2012}. Note that $s(n)=\sum_{d\div n} c(d)$
($n\in \N$).

In this paper we are concerned with the asymptotic properties of the
Dirichlet summatory functions of $s(m,n)$, $c(m,n)$, $s(n)$ and
$c(n)$. As far as we know, no such results are given in the
literature. The only existing asymptotic results for the number of
subgroups of finite Abelian groups having rank two concern another
function, namely $t_2(n)$, see Section \ref{Sect_Appendix}. We
establish asymptotic formulas for the sums $\sum_{m,n\le x} s(m,n)$,
$\sum_{m,n\le x} c(m,n)$, $S^{(2)}(x) :=\sum'_{m,n\le x} s(m,n)$,
$C^{(2)}(x) :=\sum'_{m,n\le x} c(m,n)$, $\sum_{n\le x} s(n)$ and
$\sum_{n\le x} c(n)$, where $\sum'$ means that summation is
restricted to $\gcd(m,n)>1$. Here $S^{(2)}(x)$ and $C^{(2)}(x)$
represent the number of subgroups, respectively cyclic subgroups of
the groups $\Z_m \times \Z_n$ having rank two, with $m,n\le x$. Our
main results are given in Section \ref{Results}, while their proofs
are contained in Section \ref{Proofs}.

We remark that a compact formula for the number $s_3(n)$ of
subgroups of the group $(\Z_n)^3$ and an asymptotic formula for the
sum $\sum_{n\le x} s_3(n)$ were given in \cite{HamTot2013}.


\section{Results} \label{Results}

The double Dirichlet series of the functions $s(m,n)$ and $c(m,n)$
can be represented by the Riemann zeta function, as shown in the
next result.

\begin{theorem} \label{Th_Dir_series} For every $z,w\in \C$ with $\Re z>1, \Re w>1$,
\begin{equation} \label{Dir_s}
\sum_{m,n=1}^{\infty} \frac{s(m,n)}{m^z n^w}=
\frac{\zeta^2(z)\zeta^2(w)\zeta(z+w-1)}{\zeta(z+w)} \,,
\end{equation}
\begin{equation} \label{Dir_c}
\sum_{m,n=1}^{\infty} \frac{c(m,n)}{m^z n^w}=
\frac{\zeta^2(z)\zeta^2(w)\zeta(z+w-1)}{\zeta^2(z+w)} \,,
\end{equation}
\end{theorem}

\begin{remark} \label{Rem_convo}
{\rm According to \eqref{Dir_s} and \eqref{Dir_c},
\begin{equation*} \label{connect_s_c}
\sum_{m,n=1}^{\infty} \frac{s(m,n)}{m^z n^w} = \sum_{m,n=1}^{\infty} \frac{c(m,n)}{m^z n^w}
\sum_{m,n=1}^{\infty} \frac{F(m,n)}{m^z n^w} \,,
\end{equation*}
where the function $F$ is defined by $F(m,n)=1$ for $m=n$ and
$F(m,n)=0$ for $m\ne n$ ($m,n\in \N$). Therefore (see, e.g.,
\cite{Tot2011} and \cite{Tot2013} for related properties of the
Dirichlet convolution of arithmetic functions of several variables
and of multiple Dirichlet series),
\begin{equation} \label{connect_s_c}
s(m,n)= \sum_{d\div \gcd(m,n)} c(m/d,n/d)  \quad (m,n\in \N) \,.
\end{equation}}
\end{remark}

\begin{theorem} \label{Th_asymptotics} For large real $x$ and every fixed $\varepsilon>0$,
\begin{equation} \label{s}
\sum_{m,n\le x} s(m,n) = x^2\,\br{\sum_{r=0}^3 A_r (\log x)^r} +
O\br{x^{\frac{1117}{701}+\varepsilon}} \,,
\end{equation}
\begin{equation} \label{c}
\sum_{m,n\le x} c(m,n) =x^2\,\br{\sum_{r=0}^3 B_r (\log x)^r} +
O\br{x^{\frac{1117}{701}+\varepsilon}} \,,
\end{equation}
where $1117/701\approx 1.5934$, $A_r, B_r$ ($0\le r\le 3$) are
constants,
\begin{equation*}
A_3= \frac1{3\zeta(2)}=\frac{2}{\pi^2}, \quad
A_2=\frac1{\zeta(2)}\left(3\gamma-1-\frac{\zeta'(2)}{\zeta(2)}
\right) \,,
\end{equation*}
\begin{equation*}
A_1= \frac1{\zeta(2)} \left(8\gamma^2-6\gamma-2\gamma_1+1-
2(3\gamma-1)\frac{\zeta'(2)}{\zeta(2)}+
2\left(\frac{\zeta'(2)}{\zeta(2)}\right)^2 -
\frac{\zeta''(2)}{\zeta(2)}\right) \,,
\end{equation*}
\begin{equation*}
B_3= \frac1{3\zeta^2(2)}=\frac{12}{\pi^4}, \quad
B_2=\frac1{\zeta^2(2)}\left(3\gamma-1-2\frac{\zeta'(2)}{\zeta(2)}
\right) \,,
\end{equation*}
\begin{equation*}
B_1= \frac1{\zeta^2(2)} \left(8\gamma^2-6\gamma-2\gamma_1+1-
4(3\gamma-1)\frac{\zeta'(2)}{\zeta(2)}+
6\left(\frac{\zeta'(2)}{\zeta(2)}\right)^2 -
2\frac{\zeta''(2)}{\zeta(2)}\right) \,.
\end{equation*}
\end{theorem}

\begin{remark} \label{Rem_error} {\rm Let $\Delta(x)$ denote the error term in the
Dirichlet divisor problem, i.e.,
\be \notag 
\Delta(x) := \sum_{n\le x} \tau(n) - x\log x-\c x \,, \ee and \be
\notag \theta_0 := \inf \{\theta:\ \Delta(x) = O(x^\theta)\} \,, \ee
for $x$ large. Then $O\br{x^{\frac{1117}{701}+\epsilon}}$ can be
readily replaced by $O\br{x^{\frac{3-\theta_0}{2-\theta_0}+\ep}}$.
Using the classic bound $\theta_0\ge\frac{1}{3}$, one obtains
$O\br{x^{8/5+\ep}}$. The hitherto sharpest result
$\theta_0\ge\frac{131}{416}$, which is due to M.~Huxley
\cite{Hux2003}, gives the $O$-term stated in Theorem
\ref{Th_asymptotics}.}
\end{remark}

\begin{remark} \label{Rem_const} {\rm The constants $A_0$ and $B_0$ can be constructed from the proof
below. They are quite complicated and hardly accessible to numerical
evaluation, since they involve \textit{inter alia} the infinite
series $\sum_{k=1}^\infty \tau(k)\Delta(k)k^{-2}$.}
\end{remark}

In order to formulate our result concerning the sums $S^{(2)}(x)$
and $C^{(2)}(x)$ some further notations are needed. For $K\in \N$ and $s\in \C$ let
\begin{equation} \label{def_F_K}
F_K(s) := \prod_{p^{\nu_p(K)}\,||\,K} \br{1-\eta_p(K)\, p^{-s}} \,,
\quad \text{where} \quad  \eta_p(K):= \frac{\nu_p(K)}{\nu_p(K)+1} \,,
\end{equation}
\be \label{alphas} \ba
\alpha_0(K) := F_K(1) &=
\prod_{p^{\nu_p(K)}\,||\,K} \br{1-\eta_p(K) \,p^{-1}}\,, \\
\alpha_1(K) :=F_K'(1) &= \sum_{p^*\div K}
\frac{\eta_{p^*}(K)}{p^*}\log p^* \prod_{p^{\nu_p(K)}\,||\,K,\ p\ne
p^*}\br{1- \frac{\eta_p(K)}{p}}\,, \ea
\ee
and let
\begin{equation} \label{betas}
\beta_0(K): =\tau(K)\alpha_0(K)\,,\quad \beta_1(K):
=\tau(K)\br{\alpha_0(K)(2\gamma-1)+\alpha_1(K)} \,.
\end{equation}

\begin{theorem}
\label{Th_rank_two} For large real $x$ and every fixed
$\varepsilon>0$,
\begin{equation} \label{asymp_rank_2}
S^{(2)}(x) := \sum_{\substack{m,n\le x\\ \gcd(m,n)>1}} s(m,n) =
x^2\,\br{\sum_{r=0}^3 C_r (\log x)^r} +
O\br{x^{\frac{1117}{701}+\varepsilon}} \,,
\end{equation}
\begin{equation} \label{asymp_cyclic_rank_2}
C^{(2)}(x) := \sum_{\substack{m,n\le x\\ \gcd(m,n)>1}} c(m,n) =
x^2\,\br{\sum_{r=0}^3 D_r (\log x)^r} +
O\br{x^{\frac{1117}{701}+\varepsilon}} \,,
\end{equation}
where $C_3=A_3$, $D_3=B_3$, $C_r=A_r-b_r$, $D_r=B_r-b_r$ ($0\le r\le
2$) with $A_r$ and $B_r$ ($0\le r\le 3$) defined in Theorem
\ref{Th_asymptotics} and $b_r$ ($0\le r\le 2$) given by \be \ba
\label{b_constans}  b_2&=\sum_{K=1}^\infty
\mu(K)\br{\frac{\bo}{K}}^2 =\prod_{p}
\left(1-\frac{4}{p^2}+\frac{4}{p^3}-\frac1{p^4} \right)\,,\\
b_1 &= \sum_{K=1}^\infty \frac{2\mu(K)}{K^2}\,\bo\br{\bi-\bo\log K}\,,\\
b_0 &= \sum_{K=1}^\infty \frac{\mu(K)}{K^2}\,\br{\bi-\bo\log K}^2\,,
\ea
\ee
using the notation of \eqref{alphas} and \eqref{betas}.
\end{theorem}

\begin{theorem} \label{Th_by_Walfisz} We have
\begin{equation} \label{s_asympt}
\sum_{n\le x} s(n) = \frac{5\pi^2}{24}x^2+ O\br{x\log^{8/3} x} \,.
\end{equation}
\begin{equation} \label{c_asympt}
\sum_{n\le x} c(n) = \frac{5}{4}x^2+ O\br{x\log^{5/3} x} \,.
\end{equation}
\end{theorem}


\section{Proofs} \label{Proofs}

\vskip2mm
{\bf Proof of Theorem \ref{Th_Dir_series}.} Applying the second formula
of \eqref{eval_s_m_n} we deduce for $\Re z, \Re w>1$,
\begin{equation*}
\sum_{m,n=1}^{\infty} \frac{s(m,n)}{m^z n^w}=
\sum_{d,a,b=1}^{\infty} \frac{\phi(d) \tau(a) \tau(b)}{(da)^z(db)^w}
= \sum_{d=1}^{\infty} \frac{\phi(d)}{d^{z+w}} \sum_{a=1}^{\infty}
\frac{\tau(a)}{a^z} \sum_{b=1}^{\infty} \frac{\tau(b)}{b^w} \,,
\end{equation*}
and using the familiar formulas for the latter Dirichlet series we
obtain \eqref{Dir_s}. The proof of \eqref{Dir_c}, based on the
second formula of \eqref{eval_c_m_n} is similar. \hfill $\Box$

\vskip2mm
\noindent {\bf Proof of Theorem \ref{Th_asymptotics}.}  We need the following result.

\begin{lemma} \label{Tmn}
For $m,n\in \N$ let
$$ T(m,n):=\sum_{\ell\,|\,\gcd(m,n)}
\ell\,\tau\br{\frac{m}{\ell}}\tau\br{\frac{n}{\ell}}\,, $$ i.e.,
$$ \sum_{m,n=1}^\infty \frac{T(m,n)}{m^z n^w} = \zeta^2(z) \zeta^2(w)
\zeta(z+w-1) $$ for $\Re z, \Re w >1$. Then for
an arbitrary fixed $\varepsilon>0$,
$$ S(x):=\sum_{m,n\le x} T(m,n) = x^2\,\br{\sum_{r=0}^3 c_r (\log x)^r} +
O\br{x^{\frac{1117}{701}+\varepsilon}}\,. $$ Here $c_3=\frac{1}{3}$,
$c_2=3\gamma-1$, $c_1=8\gamma^2-6\gamma-2\gamma_1+1$.
\label{Remark1} The constant $c_0$ can be constructed from the proof
below, but is not accessible to numerical evaluation for the reason
described in Remark \ref{Rem_const}.
\end{lemma}

\vskip2mm \noindent
{\bf Proof of Lemma \ref{Tmn}.} For $x$ large, let $1<y<x$ be a positive
real parameter at our disposal, and put $z:=\frac{x}{y}$. Further,
write $M:=\max(j,k)$ for short. Then,
\be \label{Sx} \ba S(x) =
\sum_{\substack{\ell M\le x\\ \ell, j, k \in \N}}
\ell\,\tau(j)\tau(k) &= \left\{ \sum_{\substack{\ell M\le x \\
\ell\le y}} + \sum_{\substack{\ell M\le x \\ M\le z}} -
\sum_{\substack{\ell M\le x \\ \ell\le y, M\le z}}\right\}
\ell\,\tau(j)\tau(k)\\ &=:S_1(x)+S_2(x)-S_3(x)\,,
\ea \ee
say. In what follows\footnote{This arrangement implies that $O(x^\theta \log
x)$ can be replaced throughout by $O(x^\theta)$, etc.}, let $\theta$
be an arbitrary fixed real greater than $\theta_0$. Then, firstly,
\be \label{Sx3} \ba
S_3(x) &= \sum_{\ell\le y}\ell \sum_{j,k\le z}
\tau(j)\tau(k) \\ &= \br{\h y^2 + O(y)} \br{z\log z + \c z +
O(z^\theta)}^2 \\ &= \h x^2 \log^2 z + \c x^2 \log z + \h \c^2 x^2
\\ &+ O\br{\frac{x^2}{y} \log^2 x} + O\br{x^{1+\theta}
y^{1-\theta}}\,.
\ea \ee

Secondly, \def\xl{\frac{x}{\ell}}
\be \label{Sx1_1} \ba
S_1(x) &= \sum_{\ell\le y} \ell
\br{\br{\xl\log\xl+\c\xl}^2 + O\br{\br{\xl}^{1+\theta}}} \\ &=
\sum_{\ell\le y} \ell \br{\xl\log\xl+\c\xl}^2 + O\br{x^{1+\theta}
y^{1-\theta}}\,.
\ea\ee

By a straightforward computation,
\be \label{Sx1_2} \ba
&\sum_{\ell\le y} \ell \br{\xl\log\xl+\c\xl}^2 \\
=&\ x^2 \sum_{\ell\le y}\frac{\log^2\ell}{\ell} -2x^2(\log x+\c)
\sum_{\ell\le y}\frac{\log\ell}{\ell}\\ &+ x^2(\log x+\c)^2
\sum_{\ell\le y}\frac{1}{\ell}\,.
\ea \ee

By Euler's summation formula, for $r=0,1,2$, \be \notag
\sum_{\ell\le y}\frac{\log^r\ell}{\ell} = \frac{\log^{r+1}y}{r+1} +
\gamma_r + O\br{\frac{\log^r y}{y}}\,. \ee

Combining this with \eqref{Sx1_2} and \eqref{Sx1_1}, we get
\be \label{Sx1} \ba
& S_1(x) =
x^2\left(\log^2 x (\log y + \gamma) \right.
\\ &-\log x\br{\log^2 y - 2\c \log y -2\gamma\c +2\gamma_1}  \\
& + \inv{3}\log^3 y - \c \log^2 y \\ &+ \left. \c^2 \log y +
\gamma\c^2 - 4\gamma \gamma_1 +2\gamma_1+\gamma_2\right) \\ &+
O\br{\frac{x^2}{y} \log^2 x} + O\br{x^{1+\theta} y^{1-\theta}}\,.
\ea \ee

Finally, with $M:=\max(j,k)$,
\be \label{Sx2_1}
S_2(x)=
\sum_{j,k\le z}\tau(j)\tau(k)\sum_{\ell\le\frac{x}{M}}\ell =
\sum_{j,k\le z}\tau(j)\tau(k)
\br{\frac{x^2}{2M^2}+ O\br{\frac{x}{M}}}\,.
\ee

The $O$-term here contributes overall
\be \label{Oterm}
\ll x \sum_{j,k\le z} \frac{\tau(j)\tau(k)}{\sqrt{j k}} \ll x \br{\sum_{j\le
z}\frac{\tau(j)}{\sqrt{j}}}^2\ll x z \log^2 x = \frac{x^2}{y} \log^2
x\,.
\ee

Writing $S_2^*(x)$ for the main term in \eqref{Sx2_1}, we
get
\be \label{S2_star} \ba
& S_2^*(x) = \frac{x^2}{2} \sum_{j,k\le
z} \frac{\tau(j)\tau(k)}{\max(j^2,k^2)}\\  & = x^2\sum_{j\le k\le
z}\frac{\tau(k)}{k^2}\,\tau(j) - \frac{x^2}{2} \sum_{k\le
z}\frac{\tau^2(k)}{k^2} =: x^2\br{R_1(z)-\h R_2(z)}\,.
\ea \ee

Now
\be \label{R2} R_2(z) =
\frac{\zeta^4(2)}{\zeta(4)}+O\br{\frac{\log^3 x}{z}} =
\frac{5\pi^4}{72}+O\br{\frac{\log^3 x}{z}}\,.\ee Further, \be
\label{R1} \ba R_1(z) &= \sum_{k\le z} \frac{\tau(k)}{k^2}\br{k\log
k +\c k+\Delta(k)} \\ &= \sum_{k\le z}\frac{\tau(k)}{k}(\log
k+2\gamma-1) + \sum_{k\le z} \frac{\tau(k)\Delta(k)}{k^2}\,.
\ea \ee

Here the last sum equals
$$
\sum_{k=1}^\infty
\frac{\tau(k)\Delta(k)}{k^2} + O\br{z^{\theta-1}} =: C_1 +
O\br{z^{\theta-1}}\,.
$$

Moreover, using Stieltjes integral notation,
\be \label{Stieltjes} \ba
&\sum_{k\le z}\frac{\tau(k)}{k}(\log k+2\gamma-1)\\ &= \Int_{1-}^{z+}
\frac{\log u+2\gamma-1}{u} \d\br{u\log u+ \c u+\Delta(u)} \\ & =
\Int_{1}^{z}\frac{\log u+2\gamma-1}{u}\br{\log u + 2\gamma}\d u  +
\c^2 + O(z^{\theta-1}) \\ &- \Int_1^z \frac{\d}{\d u}\br{\frac{\log
u +2\gamma-1}{u}} \Delta(u)\d u\\ & = \inv{3}\log^3
z+\br{2\gamma-\h}\log^2 z + 2\gamma\c\log z + C_2 +
O\br{z^{\theta-1}}\,,
\ea \ee
where
$$ C_2 := \c^2 + \Int_1^\infty
\frac{\log u+2(\gamma-1)}{u^2} \Delta(u) \d u\,.
$$

Putting together \eqref{Sx2_1} - \eqref{Stieltjes}, and recalling that
$z=\frac{x}{y}$, we arrive at
\be \label{Sx2} \ba
S_2(x) &= x^2\left(\inv{3}\log^3 z+\br{2\gamma-\h}\log^2 z + 2\gamma\c\log z
\right. \\ & +\left. C_1+C_2-\frac{5\pi^4}{144}\right) +
O\br{\frac{x^2}{y} \log^2 x} + O\br{x^{1+\theta} y^{1-\theta}}\,.
\ea \ee

Finally, using \eqref{Sx3}, \eqref{Sx1}, and \eqref{Sx2} in
\eqref{Sx}, an involved but straightforward calculation yields
$$
S(x) = x^2\,\br{\sum_{r=0}^3 c_r (\log x)^r} +
O\br{\frac{x^2}{y} \log^2 x} + O\br{x^{1+\theta} y^{1-\theta}}\,,
$$
with $c_1,c_2,c_3$ as stated in Lemma \ref{Tmn}. Balancing the two
$O$-terms here, the optimal choice is
$y=x^{\frac{1-\theta}{2-\theta}}$. This completes the proof of Lemma
\ref{Tmn}. \hfill $\Box$

Now use that (cf. Theorem \ref{Th_Dir_series} and Remark \ref{Rem_convo}),
\begin{equation*}
s(m,n) = \sum_{d\div \gcd(m,n)} \mu(d) T(m/d,n/d) \quad (m,n\in \N) \,.
\end{equation*}

We deduce
\begin{equation*} \ba
\sum_{m,n\le x} s(m,n) & = \sum_{d\le x} \mu(d) \sum_{a,b\le x/d}
T(a,b) \\ & = \sum_{d\le x} \mu(d) \left(\left(\frac{x}{d}\right)^2
\sum_{r=0}^3 c_r \left(\log \frac{x}{d}\right)^r +
O\left(\left(\frac{x}{d}\right)^{\frac{1117}{701}+ \varepsilon}\right) \right)\\
& = x^2 V(x)+ O\left(\sum_{d\le x}
\left(\frac{x}{d}\right)^{\frac{1117}{701}+\varepsilon}\right) \,,
\ea
\end{equation*}
where the error term is $O\left(x^{\frac{1117}{701}+\varepsilon}
\right)$ and
\begin{equation*}
V(x)= \left(c_3 \log^3 x + c_2 \log^2 x +c_1\log x +c_0\right)
\sum_{d\le x} \frac{\mu(d)}{d^2}
\end{equation*}
\begin{equation*}
- \left(3c_3 \log^2 x + 2c_2 \log x + c_1\right) \sum_{d\le x}
\frac{\mu(d)\log d}{d^2}
\end{equation*}
\begin{equation*}
+ \left(3c_3 \log x +  c_2 \right)
\sum_{d\le x} \frac{\mu(d)\log^2 d}{d^2} -c_3\sum_{d\le x}
\frac{\mu(d)\log^3 d}{d^2}
\end{equation*}
\begin{equation*}
= \left(c_3 \log^3 x + c_2 \log^2 x +c_1\log x+ c_0 \right)
\left(\frac1{\zeta(2)} + O\left(\frac1{x}\right) \right)
\end{equation*}
\begin{equation*}
- \left(3c_3 \log^2 x + 2c_2 \log x + c_1\right) \left(
\frac{\zeta'(2)}{\zeta^2(2)} +O\left(\frac{\log x}{x}\right) \right)
\end{equation*}
\begin{equation*}
+ \left(3c_3 \log x +  c_2 \right) \left(
\frac{2(\zeta'(2))^2-\zeta''(2)\zeta(2)}{\zeta^3(2)}+
O\left(\frac{\log^3 x}{x}\right)\right)
\end{equation*}
\begin{equation*}
-c_3 \left( c^* + O\left(\frac{\log^3 x}{x}\right) \right) \,,
\end{equation*}
with a certain constant $c^*$, leading to the asymptotic formula
\eqref{s}.

From \eqref{connect_s_c} we deduce by M\"obius inversion that
\begin{equation*}
c(m,n)=\sum_{d\div \gcd(m,n)} \mu(d)s(m/d,n/d) \,,
\end{equation*}
and obtain
\begin{equation*}
\sum_{m,n\le x} c(m,n)= \sum_{d\le x} \mu(d) \sum_{a,b\le x/d}
s(a,b) \,.
\end{equation*}

Applying now the formula \eqref{s}, similar computations show the
validity of \eqref{c}. \hfill $\Box$

\vskip2mm \noindent
{\bf Proof of Theorem \ref{Th_rank_two}.}
Obviously, by \eqref{eval_s_m_n},
\begin{equation} \label{gcd_m_n}
\sum_{\substack{m,n\le x\\ \gcd(m,n)>1}} s(m,n) = \sum_{m,n\le x}
s(m,n) - \sum_{\substack{m,n\le x\\ \gcd(m,n)=1}} \tau(m)\tau(n) \,.
\end{equation}

In order to find an asymptotics for the last sum we need an auxiliary
result.

\begin{lemma}\footnote{It is possible or even likely,
that this result or even a sharper assertion is contained in the
literature. However, the authors' attempts to find it were not
successful.} \label{Umn} For $K$ a positive integer and $Y$ a large
real variable satisfying $K\le Y^9$, it holds true that
$$
\sum_{n\le Y} \tau(K n) = \beta_0(K) Y\log Y + \beta_1(K) Y +
O\br{Y^{1/3+\ep}}\,,
$$
for every fixed $\ep>0$, uniformly in $K$, with the notations
given by \eqref{alphas} and \eqref{betas}. Note that $ \beta_0(K)\,,\ \beta_1(K) \ll
K^{\varepsilon}$.
\end{lemma}

\noindent {\bf Proof of Lemma \ref{Umn}.} We start from the
formula\footnote{The authors are grateful to Professor A.~Ivi\'c for
directing their attention to this identity.} (cf.~E.~C.~Titchmarsh
\cite[(1.4.2)]{Tit1986}) \be \notag \ba
\inv{\tau(K)}&\sum_{n=1}^\infty \frac{\tau(K n)}{n^s} =
\zeta^2(s)\,F_K(s)\, \qquad (\Re s>1)\,, \ea\ee where $F_K(s)$ is
defined in \eqref{def_F_K}. We can follow the classic example of the
deduction of \cite[Theorem 12.2]{Tit1986}, sketching only the
necessary changes. With $a_n:=\frac{\tau(Kn)}{\tau(K)} \le \tau(n)$
and \hbox{$T:=Y^{2/3}$,} Perron's formula gives \be \notag
\sum_{n\le Y} a_n = \inv{2\pi i} \Int_{1+\delta-iT}^{1+\delta+iT}
\zeta^2(s)\,F_K(s) \frac{Y^s}{s}\d s + O\br{Y^{1/3+\ep}}\,, \ee with
arbitrarily small fixed $\delta>0$. The line of integration is now
shifted to $s=-\delta+it$, $-T\le t \le T$. On the horizontal line
segments $-\delta\le\sigma\le1+\delta$, $t=\pm T$,
$$
F_K(\sigma\pm iT) \ll \prod_{p \,|\,K} p^{\delta}\le K^\delta \le
Y^{9\delta}\,,
$$
hence this brings in only a harmless factor. The residue of the integrand at $s=1$ equals
$$
\alpha_0(K) Y \log Y +\br{\alpha_0(K)(2\gamma-1)+\alpha_1(K)}Y,
$$
where $\alpha_0(K)$ and $\alpha_1(K)$ are defined by \eqref{alphas}.

Furthermore, the residue of the integrand at $s=0$ is
$$
\zeta^2(0)F_K(0)= \inv{4}\prod_{p^{\nu_p(K)}\,||\,K} \br{1-\eta_{p,K}}\
\ll 1
$$
uniformly in $K$. It remains to estimate the integral
$$
\Int_{-\delta-iT}^{-\delta+iT} \zeta^2(s)\,F_K(s) \frac{Y^s}{s}\d s \,.
$$

Expanding the product which defined $F_K(s)$ we obtain
$2^{\omega(K)}\ll K^{\delta'}$ terms (with $\omega(K)$ denoting the
number of prime divisors of $K$) of the form $B^{-s}$ where $B$ is
the product of some or all of the primes which divide $K$.

Now the proof of \cite[Theorem 12.2]{Tit1986}, which involves
$\zeta^2(s)$ alone, ultimately leads to estimates of the type
$$
\Int_1^T G(t) e^{i\Phi(t)}\d t \ll \max_{[1,T]}\abs{G(t)}\,\max_{[1,T]}\abs{\Phi''(t)}^{-1/2}\,,
$$
with $G(t), \Phi(t)$ real functions. Writing
$B^{\delta-it}=B^{\delta} e^{-it\log B}$, we see that $B^\delta\le
K^\delta$ contributes only a harmless factor, while $-t\log B$ does
not contribute at all to $\Phi''(t)$. This completes the proof of
Lemma \ref{Umn}. \hfill $\Box$

\begin{lemma} \label{Ux} For large real $x$, let
$$
U(x):=\sum_{\substack{m,n\le x\\ \gcd(m,n)=1}} \tau(m) \tau(n)\,.
$$

Then it follows that
$$
U(x) = x^2\br{b_2\log^2x+b_1\log x+b_0}+O\br{x^{4/3+\ep}}
$$
for every $\ep>0$. Here the constants $b_0,b_1$ and $b_2$ are given by \eqref{b_constans}.
\end{lemma}

\vskip2mm \noindent
{\bf Proof of Lemma \ref{Ux}.} By a familiar device usually attributed to
Vinogradov,
\be \label{vino} \ba
U(x)&=\sum_{K\le x} \mu(K) \sum_{\substack{m,n\le x\\ K\div\gcd(m,n)}}\tau(m) \tau(n)\\
&= \sum_{K\le x}\mu(K) \br{\sum_{n\le\frac{x}{K}}\tau(K n)}^2\,.
\ea\ee

The contribution of the $K$ with $x^{9/10}<K\le x$ is small: Using that
\hbox{$\tau(K n)\le\tau(K)\tau(n)$,} we get
$$
\sum_{x^{9/10}<K\le x}\mu(K) \br{\sum_{n\le\frac{x}{K}}\tau(K n)}^2
\ll \sum_{x^{9/10}<K\le x} \tau^2(K)\,\br{\frac{x}{K}}^{2+\ep}\ll
x^{11/10+\ep}\,.
$$

But if $K\le x^{9/10}$, then $K\le\br{\frac{x}{K}}^9$, thus we may apply Lemma \ref{Umn} to the
inner sum in \eqref{vino}. In this way,
\be \notag \ba
U(x) =& \sum_{K\le x^{9/10}} \mu(K)\br{\beta_0(K)\,\frac{x}{K} \log
\frac{x}{K} + \beta_1(K) \, \frac{x}{K} +
O\br{\br{\frac{x}{K}}^{1/3+\ep/2}}}^2\\ & + O\br{x^{11/10+\ep}}\,.
\ea\ee

The $O$-term here contributes overall at most
$$
\ll x^{4/3+\ep}\sum_{K\le x^{9/10}} K^{-4/3} \ll x^{4/3+\ep}\,.
$$

We claim that this implies that
\be \label{claim}\ba
U(x)=& x^2\br{\sum_{K=1}^\infty \frac{\mu(K)}{K^2}\br{\bo\log x+\bi-\bo\log
K}^2}\\ &+O\br{x^{4/3+\ep}}\,.
\ea\ee

But this is easy to see, since
\be \notag \ba
x^2& \br{\sum_{K>x^{9/10}} \frac{\mu(K)}{K^2}\br{\bo\log x+\bi-\bo\log K}^2}\\ &\ll
x^{2+\ep}\sum_{K>x^{9/10}} K^{-2+\ep}\ll x^{11/10+2\ep}\,.
\ea\ee

By an obvious calculation, \eqref{claim} completes the proof of Lemma
\ref{Ux}. \hfill $\Box$

Now Lemma \ref{Ux} and \eqref{gcd_m_n}, together with \eqref{s} and
\eqref{c}, give the asymptotics \eqref{asymp_rank_2} and
\eqref{asymp_cyclic_rank_2}, respectively.  \hfill $\Box$

\vskip2mm
\noindent{\bf Proof of Theorem \ref{Th_by_Walfisz}.} The asymptotic formulas
\eqref{s_asympt} and \eqref{c_asympt} are direct consequences of
\eqref{sss} and \eqref{ccc}, respectively, and of the known estimate
\begin{equation*}
\sum_{n\le x} \psi(n) = \frac{15}{2\pi^2}x^2+ {\cal O}(x\log^{2/3}
x)
\end{equation*}
of A. Walfisz \cite[p.\ 100]{Wal1963}. \hfill $\Box$


\section{Appendix} \label{Sect_Appendix}

Let $t_2(n)$ denote the sum of the numbers of subgroups of Abelian
groups of order $n$ having rank $\le 2$ (up to isomorphisms). Then
one has the following Dirichlet series representation, due to
G.~Bhowmik: For $\Re z>1$,
\begin{equation} \label{Dir_t_2}
\sum_{n=1}^{\infty} \frac{t_2(n)}{n^z} =
\zeta^2(z)\zeta^2(2z)\zeta(2z-1)  \prod_p \left(1+ \frac1{p^{2z}}
-\frac2{p^{3z}}\right)\,.
\end{equation}

The universality of this Dirichlet series was proved in
\cite{Lau2001}. For asymptotic properties of $t_2(n)$, based on
formula \eqref{Dir_t_2} see the papers \cite{BhoMen1997,Ivi1997}.
More generally, for finite Abelian groups $A$ let $t_r(n)=\sum_{\#
A=n, \rank(A)=r} s(A)$, where $\#A$ is the order of $A$, $\rank(A)$
is its rank and $s(A)$ stands for the number of subgroups of $G$.
See \cite{BhoWu2001} for properties of the function $t_r(n)$.

Here we give a short direct proof for \eqref{Dir_t_2}. By the
Busche-Ramanujan identity for the divisor function $\tau$ (cf. \cite{Tot2013}), the second
formula of \eqref{eval_s_m_n} can be written as
\begin{equation} \label{another_repres_s_mn}
s(m,n)= \sum_{d\div \gcd(m,n)}  d\, \tau(mn/d^2)\,,
\end{equation}
see \cite{HHTW2012}. Now, according to the definition of $t_2(n)$
and using \eqref{another_repres_s_mn},
\begin{equation*}
t_2(n)= \sum_{\substack{k\ell=n\\ k\div \ell}} s(k,\ell)=
\sum_{\substack{k\ell=n\\ k\div \ell}} \sum_{d\div k}
d\tau(k\ell/d^2)= \sum_{d^2a^2j=n} d\tau(a^2j) \,,
\end{equation*}
that is
\begin{equation} \label{ident_t_2}
t_2(n)= \sum_{d^2k=n}  d\tau(k) \tau(1,2;k) \quad (n\in \N) \,,
\end{equation}
where, as usual, $\tau(1,2;k)= \sum_{a^2b=k} 1$ (sequence \cite[item A046951]{OEIS}).
Note that $\tau(1,2;k)=\prod_p \left(\lfloor \nu_p(k)/2 \rfloor+1 \right)$. It turns out that
the function $t_2(n)$ is multiplicative and \eqref{ident_t_2} quickly leads to
the formula \eqref{Dir_t_2}.

Let $u_2(n)$ denote the total number of cyclic subgroups of Abelian
groups of order $n$ having rank $\le 2$ (up to isomorphisms), not investigated in the literature.
It follows at once from \eqref{connect_s_c} and the above results for $t_2(n)$ that
\begin{equation*}
\sum_{n=1}^{\infty} \frac{u_2(n)}{n^z} =
\zeta^2(z)\zeta(2z)\zeta(2z-1)  \prod_p \left(1+ \frac1{p^{2z}}
-\frac2{p^{3z}}\right),
\end{equation*}
and
\begin{equation*}
u_2(n)= \sum_{a^2b=n} \mu(a) t_2(b) \quad (n\in \N) \,.
\end{equation*}

\section{Acknowledgement} L.~T\'oth gratefully acknowledges support from the Austrian Science
Fund (FWF) under the project Nr. M1376-N18.


\medskip

\noindent W.~G.~Nowak \\
Institute of Mathematics, Universit\"at f\"ur Bodenkultur \\
Gregor Mendel-Stra{\ss}e 33, A-1180 Vienna, Austria \\ E-mail:
nowak@boku.ac.at

\medskip

\noindent L. T\'oth \\
Institute of Mathematics, Universit\"at f\"ur Bodenkultur \\
Gregor Mendel-Stra{\ss}e 33, A-1180 Vienna, Austria \\ and \\
Department of Mathematics, University of P\'ecs \\ Ifj\'us\'ag u. 6,
H-7624 P\'ecs, Hungary \\ E-mail: ltoth@gamma.ttk.pte.hu

\end{document}